\documentclass{article}

\usepackage{floatrow}
\newfloatcommand{capbtabbox}{table}[][\FBwidth]

\usepackage{multirow}
\usepackage{PRIMEarxiv}

\usepackage{microtype}

\usepackage{graphicx}
\graphicspath{{images/}}
\DeclareGraphicsExtensions{.pdf,.png,.jpg}

\usepackage[T1]{fontenc}    
\usepackage{url}            
\usepackage{booktabs}       
\usepackage{amsfonts}       
\usepackage{nicefrac}       
\usepackage{microtype}      
\usepackage{lipsum}
\usepackage{fancyhdr}       
\usepackage{graphicx}       
\usepackage{wrapfig}
\graphicspath{{media/}}     

\title{\LaTeX}  

\title{Investigating the probability of a cylindrical coin landing on its side}

\author{
  Anton Gaek\\
  Oxford International College\\
  Oxford, England\\
  \texttt{tonygaek@gmail.com} \\
   \And
  Artem Sukhov\\
  Moscow Institute of Physics and Technology\\
  Dolgoprudny, Russia\\
  \texttt{sukhov.ad@phystech.edu} \\
}

\begin{document}

\maketitle

\begin{abstract}
The problem of creating a three-sided dice with the probability of it landing on each of its sides being equal to 1/3 has been around for many years. Various approaches have been attempted, but as different authors achieved at different results, no uniform answer has been found. In this paper, the probability of a cylinder-shaped coin landing on its side is investigated, aiming to predict the final position of the coin, based on the starting physical and geometrical parameters of the coin. Pure physical methods and computer modeling have been used to find the parameters of a "fair" coin and to determine which parameters of the system affect the probability of the coin landing on an edge, and in what way. Statistical analysis of the situation has been made, to ensure that the probabilities we are trying to find, correspond to the appropriate physical events. Overall, two main possibilities for modeling the situation are considered: single-axis rotation (flat rotation) and the superposition of rotation against multiple axes (volumetric rotation). For each, a model was built and it was compared with computer modeling and real-life experiments. Dynamical and energetic approaches have been considered, and a model of the coin during its impact with the ground has been created. A robot that is capable of launching cylindrical coins, picking them up, and recording the result of the experiment, controlling starting parameters precisely, and distributing them to our needs, and automatically record the results of the experiment, was created. In conclusion, the parameters of a "fair" coin are provided, which have even probabilities landing on all its sides, as well as the dependencies of the probability of a coin landing on its edge, on different parameters. 
\end{abstract}

\keywords{Probability \and Impact model  \and Robotics \and Coin toss}

\section{Introduction}
Landing a coin on its side is often associated with the idea of a rare occurrence. What should be the physical and geometrical characteristics of a cylindrical dice so that it has the same probability of landing on its side and one of its faces? In this work we are going to investigate relevant parameters that affect the probability of our coin landing on its side; work out those probabilities depending on the parameters of the coin and starting parameters of an experiment; consider different approaches to the problem; create a model of this coin’s movement in the space ad visualize this theoretical model with help of a computer; compare results of modeling, theoretical model, and experiment.

Probabilities in physical processes and so-called chance devices have been studied for a long time. Mechanical devices that are controlled by laws of physics such as dice, roulette wheels, and others have been studied by classical authors like Poincaré [Poincare, 1912; Mazliak, 2013], at around 1900 who explored how some objects obeying the laws of physics can act randomly, while seems that its behavior is completely deterministic. One of the main conclusions he arrived at, was that the randomness of the process emerges from the smallest deviations in starting conditions of the system being analyzed and that even the tiniest percentage errors in those conditions can lead to the final state of the system being completely different from what was predicted. Similar ideas have been presented in the works of another classic - Smoluchowski [Smoluchowski, Natunvissenschaften, 1918]. Almost later, in 1986, Keller [Keller, 1986] in his work on roulette wheels, and more importantly for this paper, coin tosses, presented his answers to the question of the randomness of seemingly predefined physical events and how to calculate the probabilities of these random events. By considering the mechanics of a coin toss, more specifically the differential form of Newton's equations of motion, and simple equations for rotation, he analyses the final equations to understand from them how the coin lands, and then finally gets a probability density function for the way the coin flipped. He gave his formula for the probability of heads and explained that the outcome is random, because of his assumption that the initial conditions are random. In his model, the probability did not exceed one-half, as he considers the angular and linear velocities tending to infinity.
In some later works, such as [Mahadevan, Yong, 2014] written in 2014, geometrical and more complex dynamical approaches are considered, taking dynamics of the spin and precession into account with a conclusion that a dynamical approach works better, validating the final formula using results of simple tabletop experiments, tossing the glued together US quarters. Other works like [Stefan, Cheche, 2016] account for the bounce of the coin, considering the coin as a rigid body, using Euler parameters (normalized quaternions) to derive the full motion equations, experiments conducted in a simple setting, recording the motion on camera Similar approach is taken in [Strzalko et al., 2008]. Other authors [Scott et al., 2012] consider a more mathematical approach, using the Bayesian model, although taking physical parameters such as mass and velocity into account, not considering any physical processes, such as a model for impact with the ground or complex equations of motion; statistically analyzing different models, based on the surface area, cross-sectional length, solid angle, the center of mass or simple bounce, trying to predict which one is the most likely – in conclusion leaning towards the bounce model, but saying that more complex physical analysis would be justified as in their experience it significantly improves the model. Another work that is worth mentioning and that had an impact on our approach is [Nagler, Richter, 2008] which investigates the randomness of the coin toss, and specifically something that we believe is valuable, the borderline between the cases where the coin toss can be predicted well, and ones where it can be considered chaotic and dependence on initial conditions is particularly high and sensitive. Another paper that presents interesting results is [Hernandez-Navarro, Piñero, 2021] analyzing the partially inelastic toss, providing the connection of a coin toss to a cuboid toss, which turns out to have a lot in common.

All the works discussed earlier arrive at different conclusions; some, like [Mahadevan, Yong, 2014] for example do not consider the rebound of the coin after hitting the table; no studies seem to have validated their model using a robot that can control the starting parameters reliably, accompanied by computer engine modeling and other experiments – all of which was done in this research to rigorously prove that the formula acquired in this paper for a three-sided coin is correct. As with any statistical problem, it is important to consider how the experiments and such are carried out, as different types of experiments may provide different kinds of error – an experiment when throwing a coin by hand has a human factor and corresponding biases, computer engines may not be an accurate representation of the real world, and their random number generators might not be truly random, robots might work exactly as they were intended to. This research fully considers different experiment settings while some authors [Mahadevan, Yong, 2014] do not conduct experiments carefully. But in the case of all these fitting the theoretical model, there is a strong ground to believe that the model is correct. Additionally, no other research before this one considered the multiple rebounds, and as this one does, it allows to work with high velocities of the initial rebound, which can impact the probability of the coin landing on heads. Another parameter that was not usually considered in previous researches, but that is present in this one is friction between the coin and the ground, which can have an impact on the final rotation and slippage. To summarize, this research focuses on considering the psychical parameters in modeling the cylindrical, three-sided coin, its movement, and the probability of it landing and coming to rest on a specific edge. In this work, the coefficient of energy loss during the rebound, the coefficient of friction between the coin and the surface it is falling onto, and the geometrical parameters of the coin (height and radius) are considered.

To start things off, it is needed to find all the parameters of the system coin-table, in order to find dependencies on them in the future. In the first stage of the investigation, Young's module, Poisson's coefficient, etc., are going to be ignored, not including the impact with the surface in our model. Therefore, the geometrical parameters of our system are coins height – H and coins radius – R. And physical parameters are the mass of the coin – m, coefficient of energy restitution – k, and friction coefficient between the coin and the surface – µ (figure 1).

\begin{wrapfigure}{r}{0.2\textwidth} 
    \centering
    \includegraphics[width=1\textwidth]{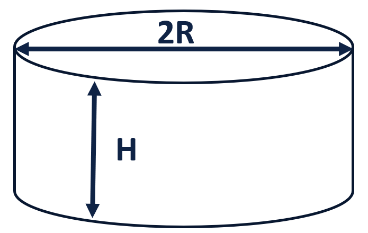}
    \caption{Coin geometrical parameters}
\end{wrapfigure}

Subsequently, let us move to the part of the problem, connected with probabilities. We are looking to express the probability of a cylinder coin landing on a particular side, from the parameters of the coin, like its dimensions, the material it is made of, etc. For the fair coin, we would like $P_{sides}$ = $P_{up}$ = $P_{down}$ = 33\% - in this case, the probabilities of landing on each "tails"\ , "heads"\ and "sides"\ would be equal. To prove our theoretical model, we are going to need to provide even distribution of starting parameters of the coin (angular velocity, angle, etc.) in our experiment. We are going to calculate the error of experimental statistics according to this formula [Reznichenko, 2017]:
\begin{equation}
\frac{\Delta P}{P}\approx \frac{1}{\sqrt{N}}\\  
\end{equation}
Another thing that we need to consider before starting to study the problem, is that we do not accept any determined situations, meaning that we can not say that everything in the coin toss is predetermined, and this is where all accidentalness of the coin toss comes from. This is exactly how our coin does its main purpose – to land on a random side. In this case, the source of uncertainty is itself the starting conditions of the coin and the range of values they can take. 

\begin{figure}
    \centering
    \includegraphics[width=0.2\textwidth]{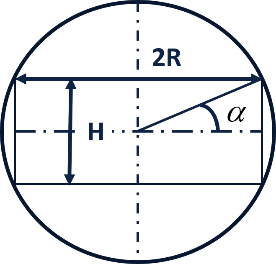}
    \caption{Flat rotation geometrical parameters}
\end{figure}

As it was mentioned earlier, the problem is going to be considered from a couple of different sides, and then decide which one is better in what instances. There are two types of rotation – flat rotation (as a coin does), and volumetric rotation (as a regular dice does). Flat rotation is applicable when a body rotates only in two different axes (2D rotation), and volumetric rotation considers rotation in all three possible directions (3D rotation).

\begin{figure}
    \centering
   \includegraphics[width=0.5\linewidth]{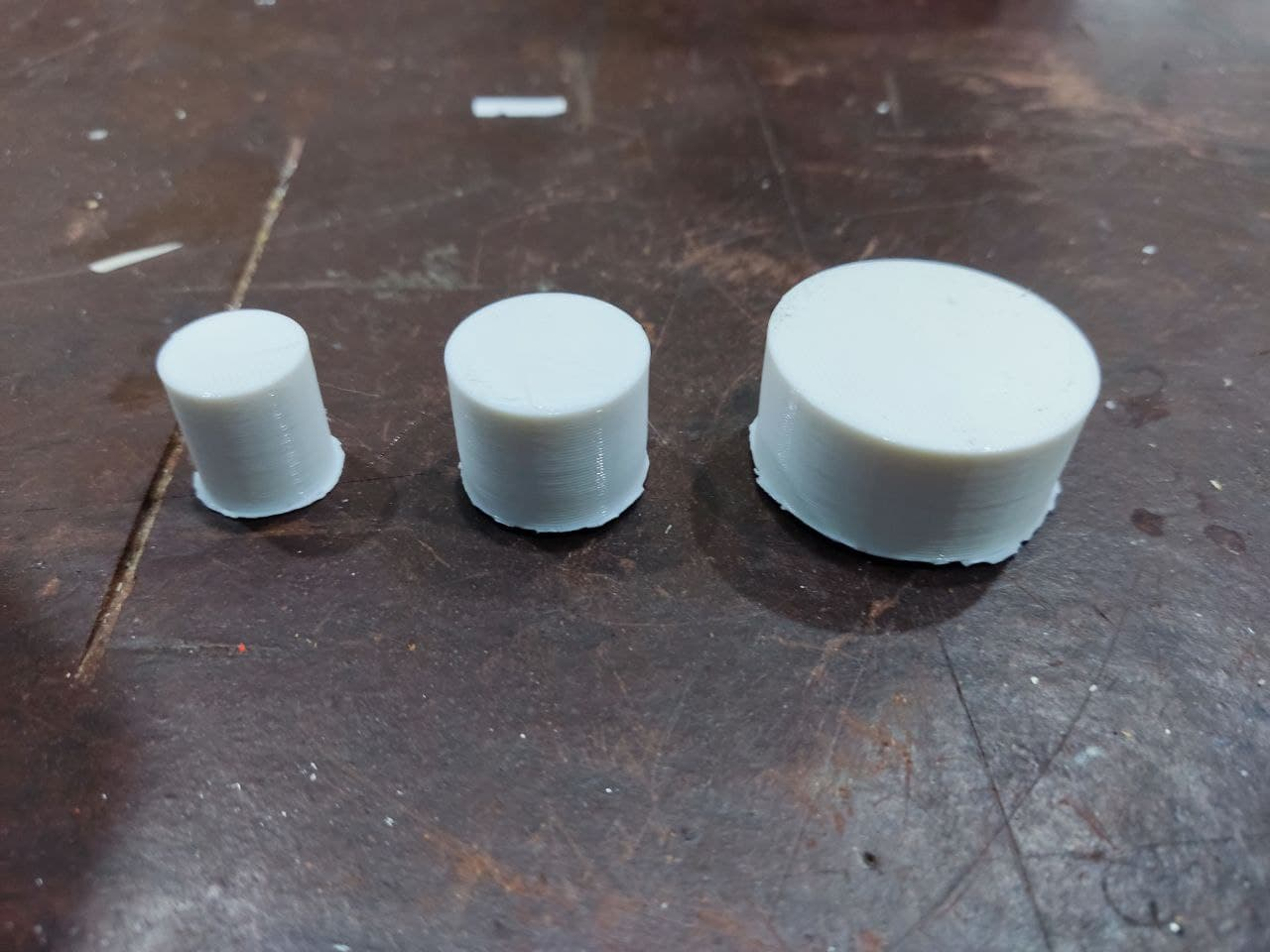}
   \caption{3D-printed coins (before the grinding process) }
\end{figure}

\section{Flat rotation (2D rotation)}
It is possible to look at the flat-2D rotation, as a rotation inscribed in a circle. Assuming an impact with the surface absolutely non-elastic, it is obvious that the probability of a coin landing on a certain side equals the probability of a circle that the rotation is inscribed in, landing on a corresponding sector (arc). From now on, we are always going to be calculating and talking about the probability of the coin landing on its side (not heads or tails). Then, from simple geometry, we can easily calculate that the probability of the coin landing on a corresponding arc is 
\begin{equation}
P=\frac{2\alpha}{\pi},
\end{equation} where $\alpha= \arctan \frac{H}{2R}$ (figure 2).


\subsection{Experimental setup}
Now, let us check this theoretical model. During the work on this problem, we created and tried different methods of making an experiment. One of them was a robot that we created, to help us make experiments. It allows us to vary the coin's angular velocity from 0 to 5 rotations per second; starting angle from 0 to 180 degrees, and overall, we get an even distribution between angles and angular velocities. The robot consists of three parts, one that throws the coin (figure 3). Then, one that drives to the coin and takes a photo of it, sends it to the computer (figure 4), so we can then use this photo to determine on what side the coin has landed using OpenCV. 
The third part of the robot elevates the coin back to the throwing part, and the cycle repeats (figure 5). The entire work cycle takes up to around 10 seconds, therefore we can make about 4000 throws per day. We have also checked if our robot actually works: we gave it a regular coin and a regular dice and made it throw them, and after about 1000 throws for both coin and dice, we got probabilities of tales about 50 \% for the coin, and probability about 16.6\% of “1” for the dice, and this shows that the robot actually works as it was intended.

\begin{figure}
    \centering
    \includegraphics[width=0.5\linewidth]{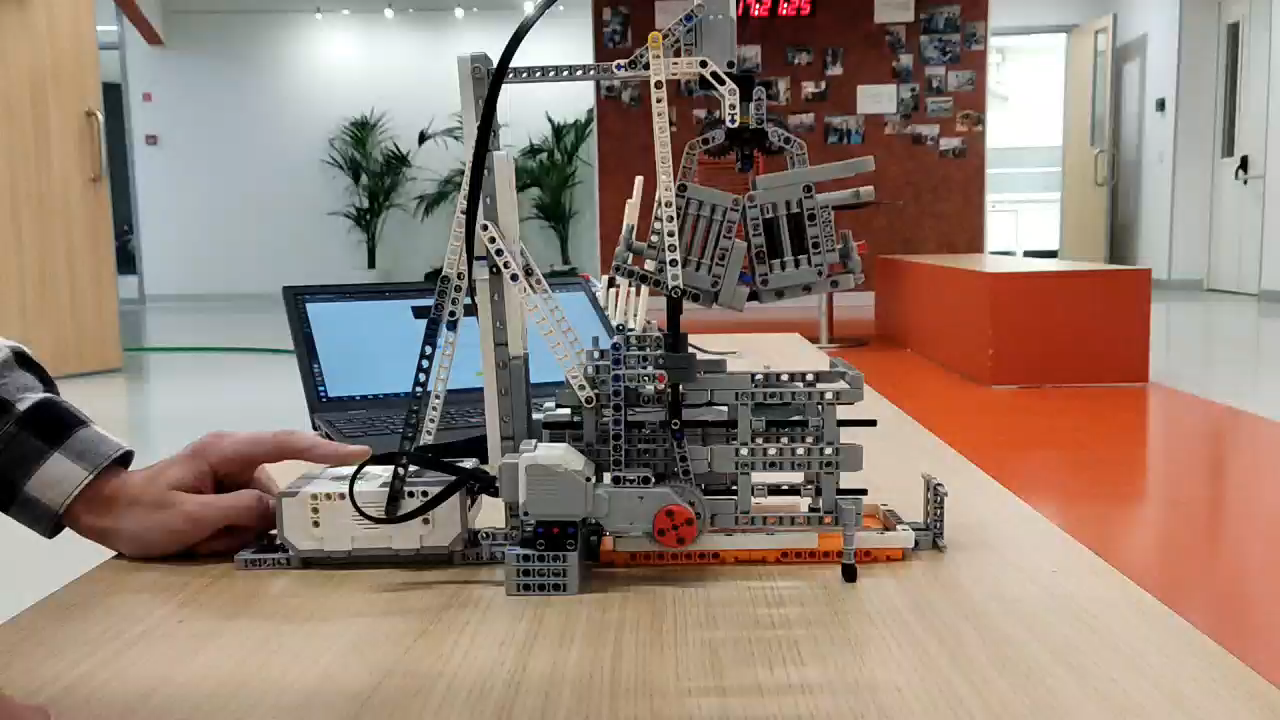}
  \caption{Throwing part of the robot}%
\end{figure}

\begin{figure}
    \centering
    \includegraphics[width=0.5\linewidth]{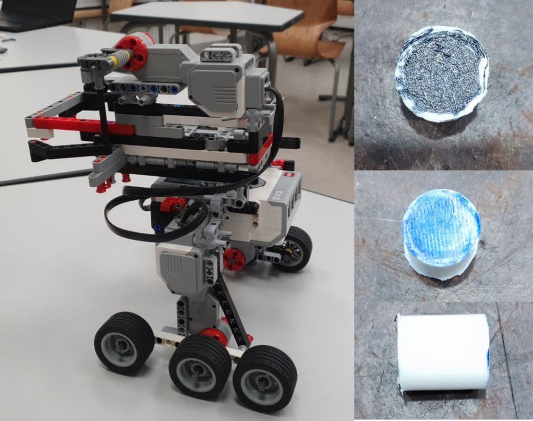}
  \caption{Photographing part of the robot (left), coins with markers for OpenCV (right)}%
\end{figure}

\begin{figure}
    \centering
    \includegraphics[width=0.5\linewidth]{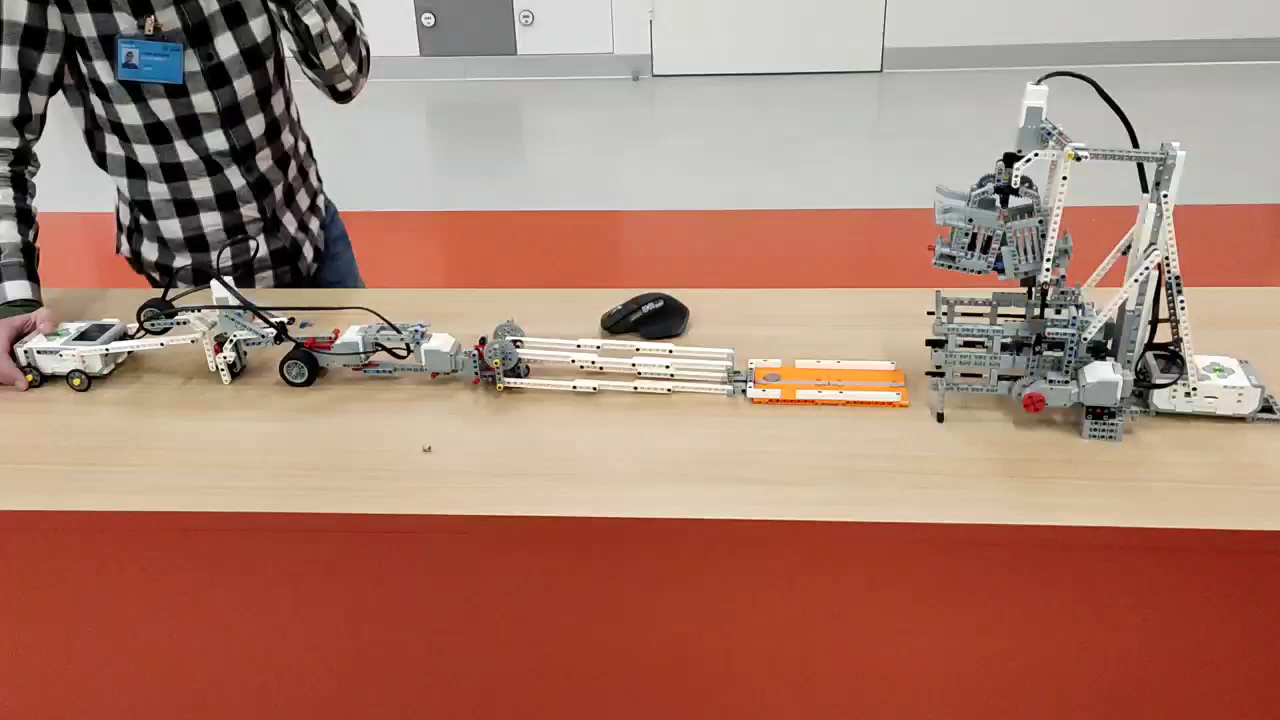}
  \caption{Lifting part of the robot to relaunch the process}%
\end{figure}

We printed out our coin on a 3D printer (figure 6), where the mistake of printing was 0,1 mm. In this way, we managed to find out the probability of our coin landing sides, from the parameters of the coin, and its mass in particular. We printed out several coins of the same volume and from, but different masses (it is easy to set up a 3D printer to print with different densities). Then we experimented with these coins with help of our robot and threw the coins by hand. All these experiments showed us that the probability of the coin landing on different sides does not depend on its mass if other parameters remain unchanged.

Now we can present the results of the first experiments, that were checking the theory of 2D rotation (2). On this graph, there is the dependence on the probability of the coin landing on its side from ratio $\frac{H}{R}$ (figure 7).

\begin{figure}[!ht]
    \centering
    \includegraphics[width=0.7\linewidth]{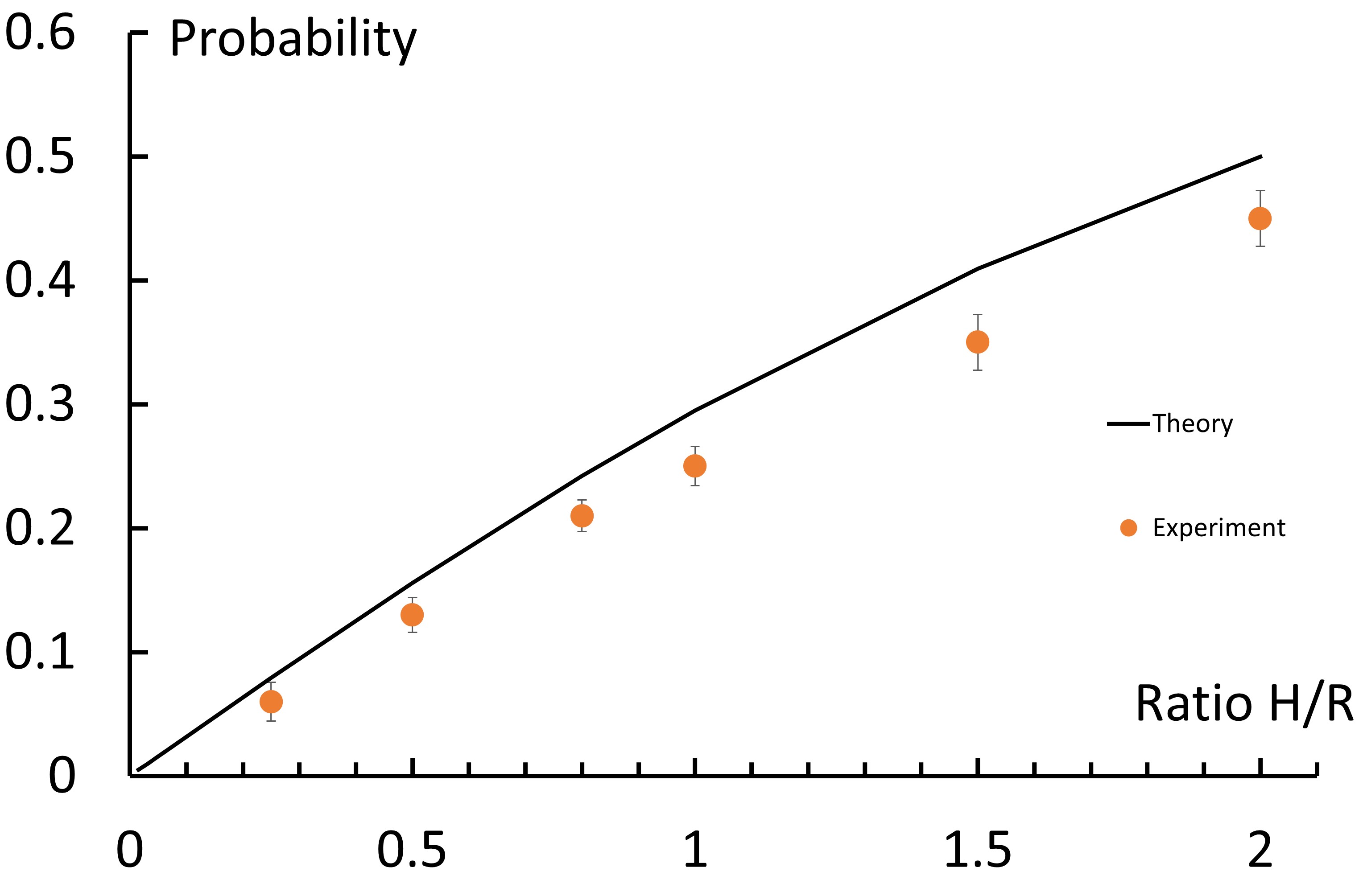}
    \caption{Probability of dropping the side of the cylinder from the H/R ratio}
\end{figure}


We can clearly see, that in the experiment there is always a systematic shift down from the theoretical curve. We need to figure out, why this happens.

Let us provide a brief summary, of everything that we have already described. We found out that the probability of the coin landing on its side does not depend on the coin's mass. Meaning that only the coefficient of momentum restitution k and coefficient of friction µ remain from the physical parameters that we wanted to include in our model. At this point, we are working with a flat model, meaning we do not want to think about friction just yet. However, to include the coefficient of momentum restitution k in our model, we need to analyze the impact of the coin with the surface.

\subsection{Impact model}

\begin{figure}[!ht]
    \centering
    \includegraphics[width=0.35\textwidth]{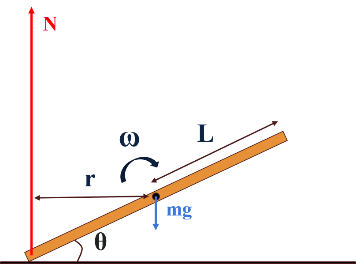}
  \caption{Coin hitting the ground left side}%
\end{figure} 

Let us consider the sum of all torques, relative to the center of mass, of a falling coin: $M=\sum_{}^{}[\overline r \times \overline F] = NL\cos \theta$, where $r$ is shown in figure 8; F under the sum are the forces that act on the coin, N is the normal reaction force shown in figure 8. In the first instance (figure 8) the coin will start spinning, and its translational energy will transfer to rotational energy. However, in the symmetric situation (figure 9) $M=\sum_{}^{}[\overline r \times \overline F]= -NL\cos \theta$, everything will happen in an opposite way, as torque of force now stops the rotation, and the rotational energy will transfer to translational energy.

\begin{figure}[!ht]
    \centering
    \includegraphics[width=0.35\textwidth]{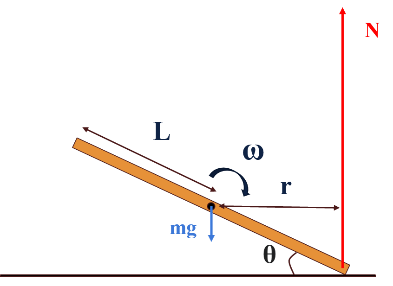}
  \caption{Coin hitting the ground right side}%
\end{figure}

Now we can use an energy conservation law for our coin, for three states of our coin: 1 – point of launch/highest elevation since the previous impact; 2 – right after the impact with the surface; 3 – point of highest elevation after the impact with the surface (figure 10)

\begin{figure}[!ht]
    \centering
    \includegraphics[width=0.35\textwidth]{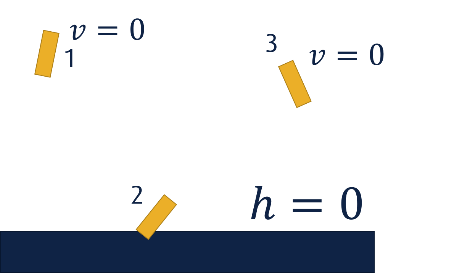}
  \caption{Stages of movement}%
\end{figure}
\begin{equation}
    k\cdot mgh_1+\frac{I\omega_1^2}{2}=\frac{mv_2^2}{2}+\frac{I\omega_2^2}{2}=mgh_3+\frac{I\omega_2^2}{2},
\end{equation}
where k – empiric coefficient of restitution; m – coins mass; $\omega_i$ – angular velocity in the i-th state of the coin (i is 1, 2 or 3); $v_i$ – velocity in the i-th state of the coin; $h_i$ – height above the surface in i-th state of the coin. From this formula, we can find $h_3$
\begin{equation}
    h_3=kh_1+\frac{I\omega_1^2 - I\omega_2^2}{2mg}.
\end{equation}
In this formula, we do not know only $\omega_2$ – the angular velocity right after the impact with the surface. We can express $\omega_2$ from $\omega_1$, so we just need to find how the angular velocity changes during an impact with the surface. We know that $I\epsilon=M$, where $\epsilon=\frac{d\omega}{dt}$. So:
\begin{equation}
    I\Delta \omega = \int\limits_0^\tau Mdt.
\end{equation}
Where $\tau$ – is impact time. Using the formula for torque M, we can simplify our expression 
\begin{equation}
    \Delta \omega=\frac{H\cos \theta  \int\limits_0^\tau Ndt}{I}.
\end{equation}

\begin{wrapfigure}{!ht}{0.2\textwidth} 
    \includegraphics[width=0.7\textwidth]{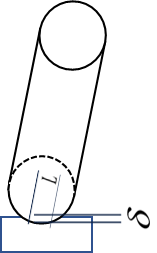}
    \caption{\footnotesize Impact model parameters}
\end{wrapfigure}

Now we only need to calculate the impact with the surface. We are going to use the results from this article [Harris, Kotzalas, 2001]. We can find these formulas in this article:
\begin{equation}
    N=K\delta ^n ; \; n = \frac{10}{9} \;\;\;\;\;\;\;\; K=\frac{1}{1.36^n}E^* L^\frac{8}{9} \;\;\;\;\;\;\;\; \frac{1}{E^*}=\frac{1- \upsilon ^2}{E},
\end{equation}

where N – impact force; K – stiffness parameter, n – Non-linearity coefficient; L – contact length; $\delta$ – relative indentation; E – Young's module; $E^*$– modulus composition; $\upsilon$ – Poisson's coefficient (figure 11).

They provide a complex model for the impact of an object with the ground, as the object barely sinks into the ground. From all these formulas, we can conclude one thing important to this research – that we can assume the force harmonic from time (graph with according parameters can be seen in figure 12). Then $N=N_{max} \sin \omega _c t;$, where $\omega _c = \frac{\pi}{\tau}$. We know that the time of the impact $\tau$ is about one millisecond. Plugging this into the formula for angular velocity change gives us
\begin{equation}
    \Delta \omega = \frac{H\cos \theta  \:  \cdot \:2N_{max} \tau }{I\pi}=\eta \frac{H \cos \theta}{I},
\end{equation}
where $\eta$ – constant material coefficient.

\begin{figure} [!ht]
    \centering
    \includegraphics[width=0.7\linewidth]{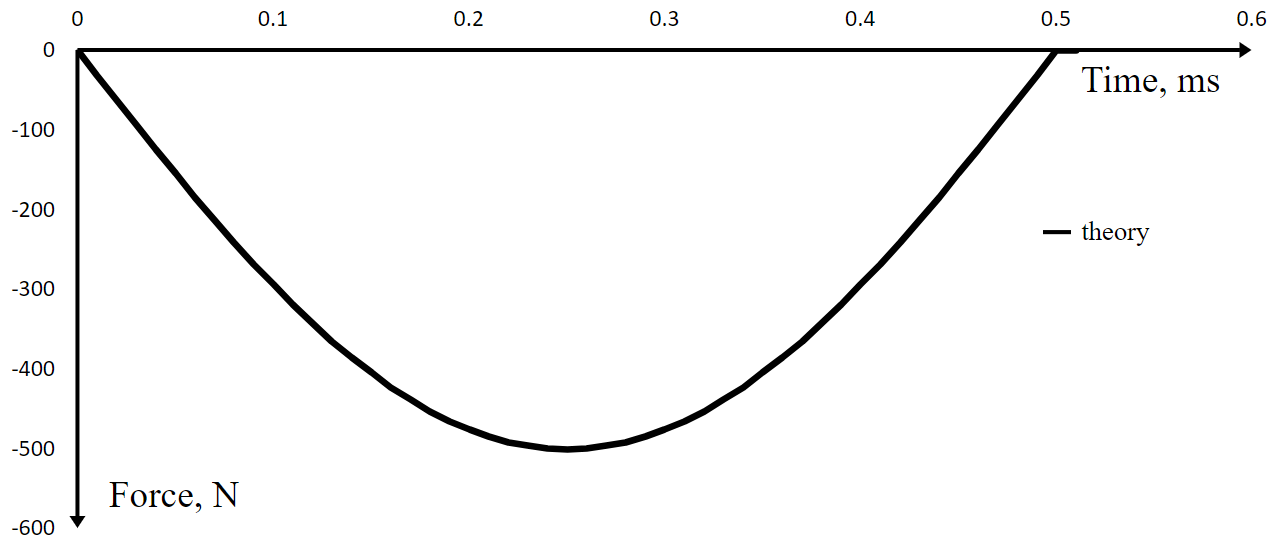}
  \caption{Harmonic representation of force}%
\end{figure}

\begin{figure}[!ht]
    \centering
    \includegraphics[width=0.5\linewidth]{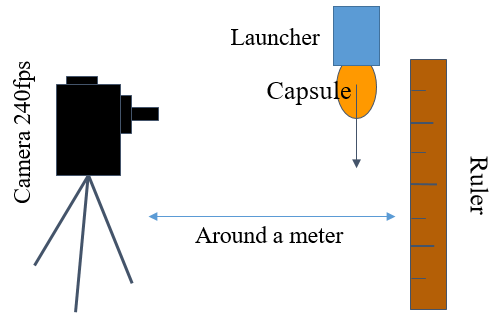}
  \caption{Experimental setup}%
\end{figure}

\subsection{Experiment for surface impact model}
For this experiment, we threw the coin down from different heights, captured the entire process on camera (240FPS), and then processed this footage in tracker [Brown, Cox, 2009]. In this way, we can check the theoretical model of the coin's movement, with experiments. You can see the scheme of the experimental setup in figure 13.

As we can see from these graphs, they fit the theory pretty well, meaning that our formula is correct and working (figure 14). Now we can make use of the formulas that we got and model the movement of the coin. Our program followed this algorithm: first of all, our code calculates the second Newton's law for translational and rotational movement, then during an impact with the surface we use the model that we described earlier to calculate new parameters of the coin (velocity, angular velocity, etc.). Then we repeat the entire process while the energy of the coin is not 100 times less than its initial energy, meaning that the coin would have stopped moving. This algorithm can help us compare theoretical models and experiments. In the end, we get the graph in figure 15.

\begin{figure}[!ht]
\centering
\includegraphics[scale=0.6]{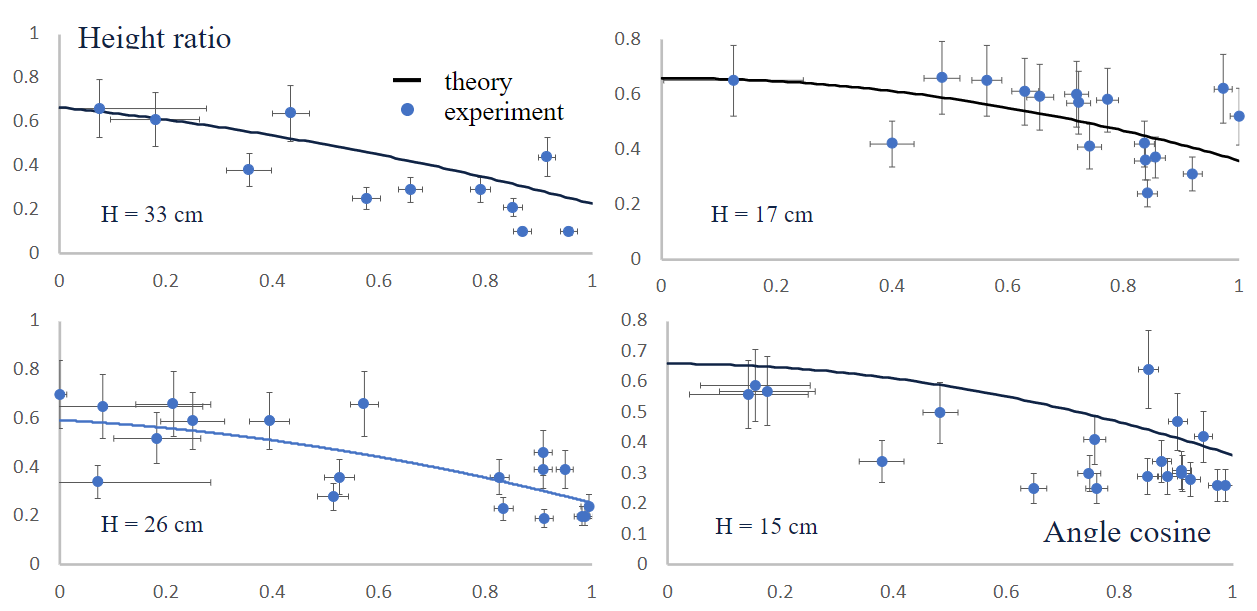}
\caption{Experiments for impact model}
\end{figure}


Let us compare the experiment with modeling. As we can see, the first experimental points fit modeling really well (orange points), and this allows us to find the ratio for fair coin $\frac{H}{R}=1.5\pm 0.05$, however, when $\frac{H}{R}$ gets bigger, experiment and modeling become different (blue points). This is because the last experimental points were made by throwing the coin by hand, meaning that there is rotation in all three axes. We could also see that the elasticity coefficient lowers the probability of the coin landing on its side.

\begin{figure}[!ht]
\centering
\includegraphics[width=0.7\linewidth]{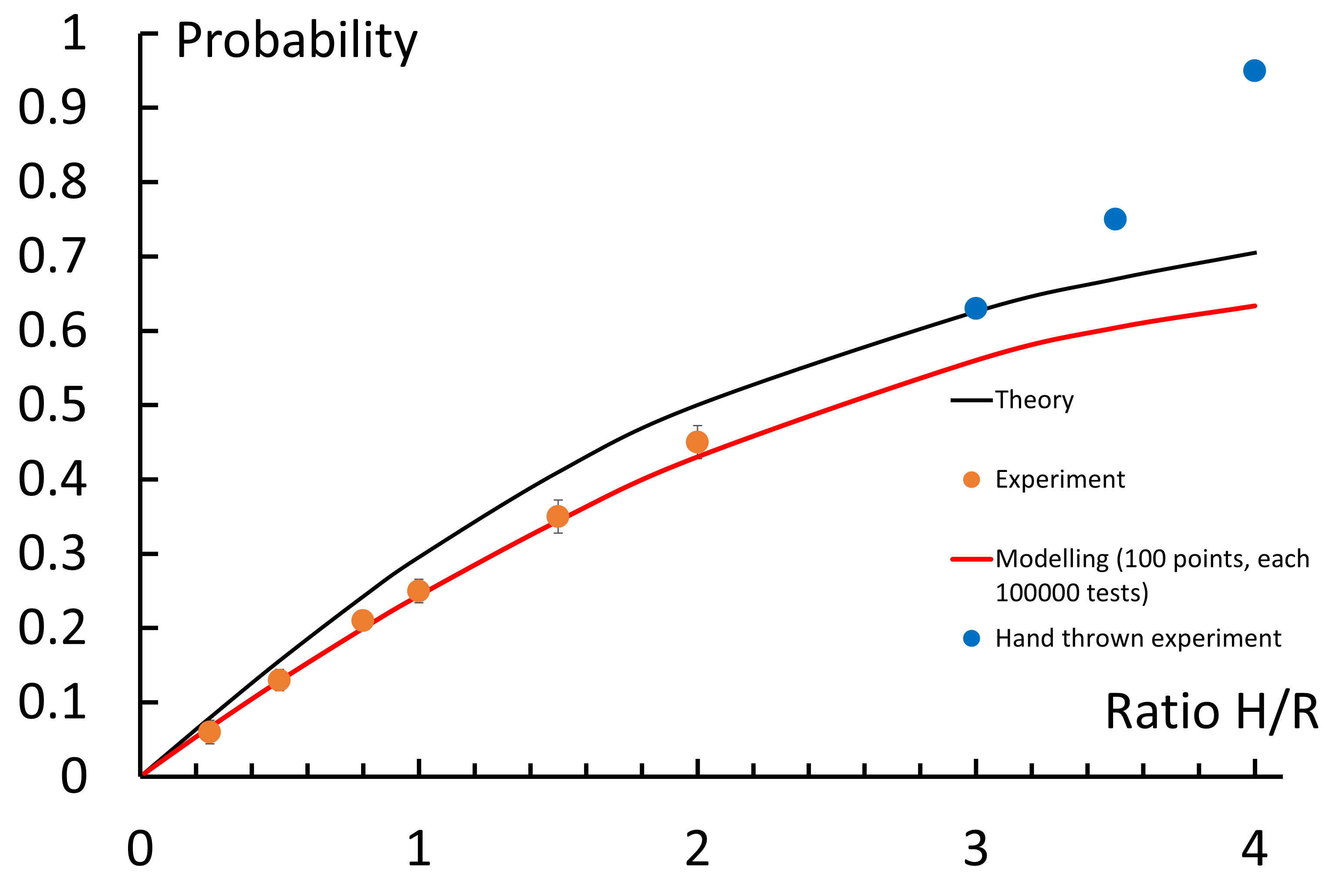}
\caption{Probability from H/R ratio graph for flat rotation}
\end{figure}


\section{Volumetric rotation (3D rotation)}
Now we can investigate the volumetric rotation model, rotation inscribed into a sphere. Analogical to flat rotation, if the impact was absolutely non-elastic, we can calculate areas of the corresponding sides to the coin, parts of the surface of the sphere that rotation is inscribed in. Then we need to find the surface area of a side surface $S_{side}=\pi H \sqrt{4R^2 + H^2}$, then the probability we want to find is $P=\frac{\frac{H}{R}}{\sqrt{\left ( \frac{H}{R} \right ) ^ 2 +4}}$ (figure 16).

\begin{figure} [!ht]
    \centering
    \includegraphics[width=0.25\textwidth]{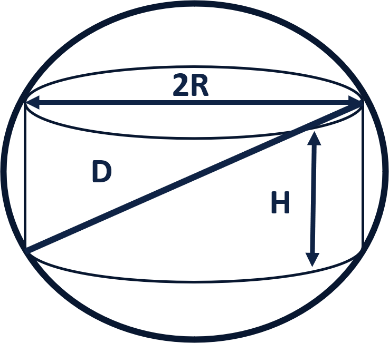}
    \caption{3D rotation parameters}
\end{figure}

We also made an experiment, throwing the coin by hand, as it guarantees randomness different from the one we can achieve using a robot or a computer simulation.

\subsection{Modelling with Unity}
To check our volumetric model theory, we have created a computer model [Sukhov and Gaek, 2023]. Our first model for planar rotation would not work here (although we used it for Unity modelling), so we are going to use Unity to model. At this point, it is important to mention that A Unity model does not fully substitute an experiment. The physical accuracy of the Unity engine is not proven, so we can not fully rely on it when drawing a conclusion if our model is working or not. However, Unity allows us to set millions of experiments, which would obviously not be possible in real life. So when we use Unity alongside a real-life experiment, to confirm our results, it helps substantially. With Unity, we can set the span of angular velocities and other physics parameters (see table 1). You can see how the modeling looks from the side in figure 17.

\begin{figure} [!ht]
\centering
\ffigbox{%
  \includegraphics[scale=0.35]{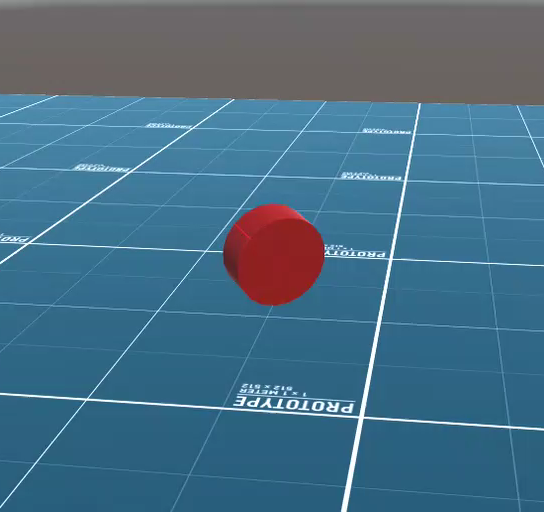}
}{%
  \caption{Coin in Unity}%
}
\capbtabbox{%
  \begin{tabular}{|p{4.5cm}||p{2cm}|}
    \hline
    \multicolumn{2}{|c|}{Unity engine parameters (for specific simulation)} \\
    \hline
    Parameter & Value\\
    \hline
    Top speed   & 20\\
    Angular speed &   50\\
    Initial height & 5 \\
    H to R ratio     & 0.8\\
    Time scale&  100\\
    N Repit& 1000\\
    Static Friction& 10\\
    Bounciness& 0.3\\
    \hline
 \end{tabular}
}{%
  \caption{Parameters of modelling}%
}

\end{figure}

In the graph, in figure 18 we can see the results of modeling in Unity, theory, and results of experiments. On this graph, each point of the experiment was acquired from at least 500 throws, that we made by hand; each point of modeling from at least 100000 throws. We throw the coin with our hand, to achieve rotation in all three axes and to provide another type of random experiment that can prove our model.  



We can see that theoretical modeling, experiment, and computer modeling all fit each other pretty well. We can find ratio $\frac{H}{R}=0.8\pm 0.05$ for a fair coin. However, there are more parameters that can influence the result.

One of those parameters is the coefficient of friction – $\mu$. We can not neglect this parameter as friction might have a crucial effect on the behavior of the coin. Using Unity, we can build a graph of the probability of a coin landing on its side from the coefficient of friction $\mu$ (figure 19).
As we can see, dependence on $\mu$ is small, but it is still there so we cannot neglect it.

\begin{figure} [!ht]
    \centering
    \includegraphics[width=0.7\linewidth]{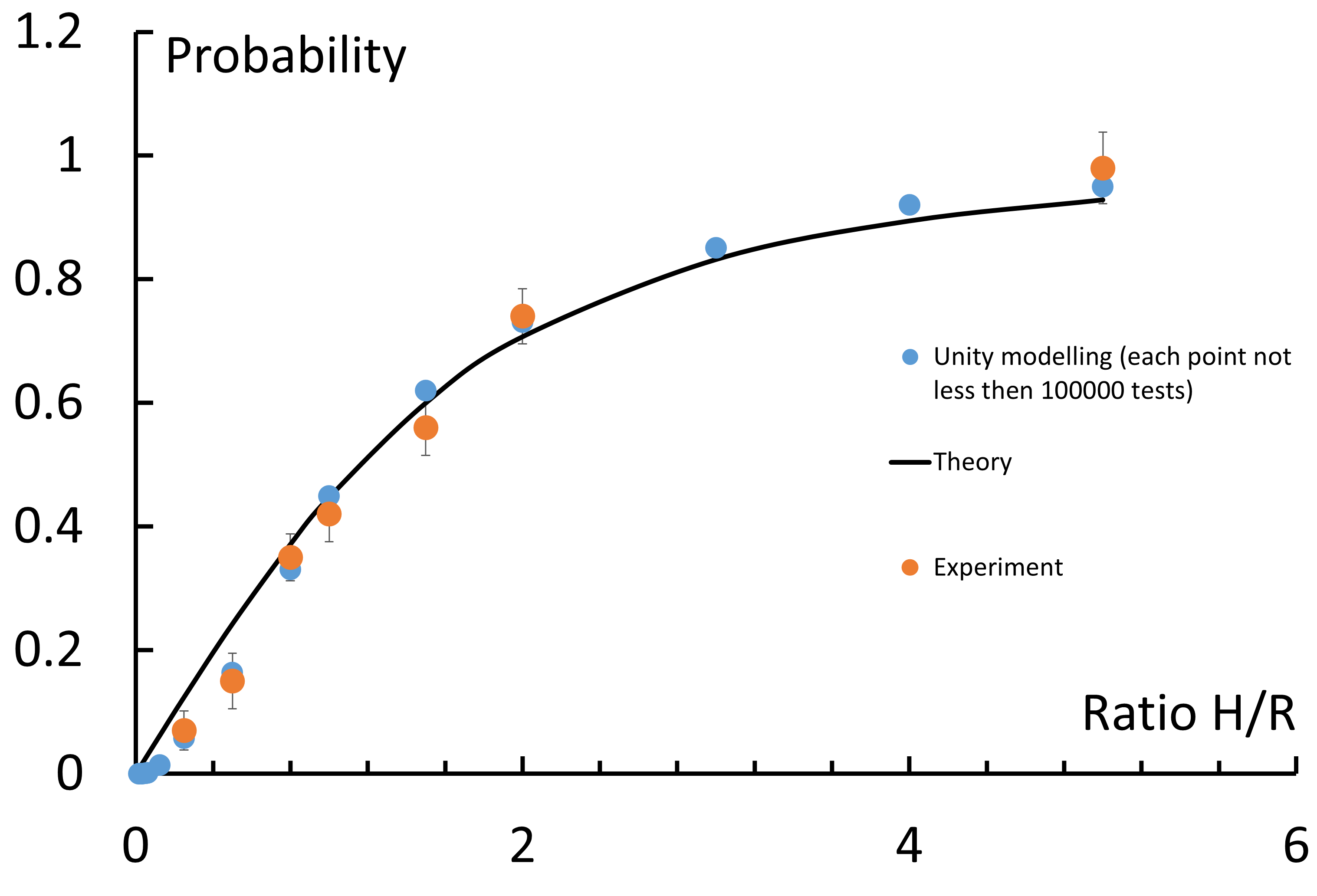}
    \caption{Probability from H/R ratio graph for volumetric rotation}
\end{figure}

\begin{figure} [!ht]
    \centering
    \includegraphics[width=0.7\linewidth]{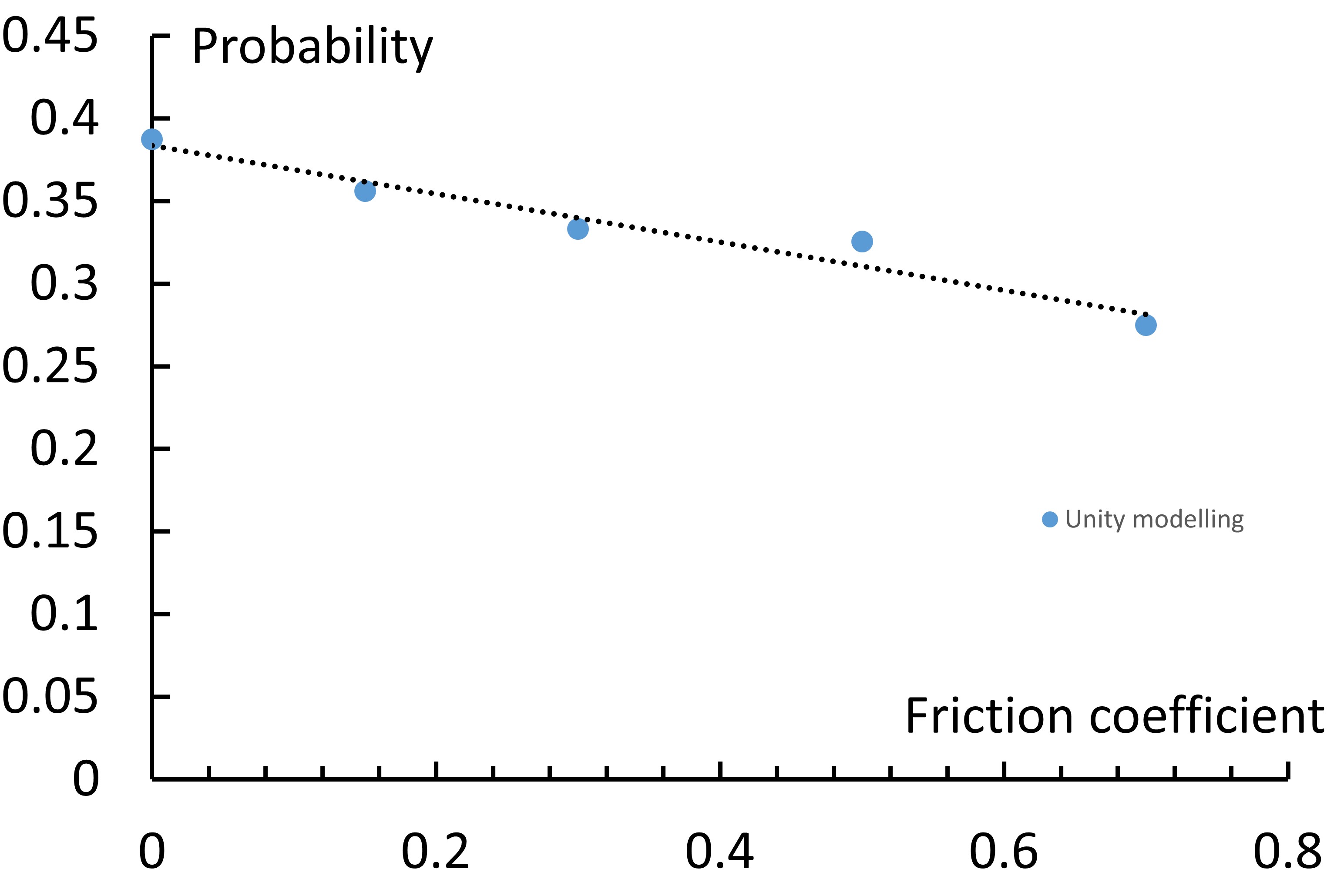}
    \caption{Probability from friction coefficient graph}
\end{figure}

We need to explain why all methods of measuring probability (all experiments, modeling, etc.), altogether prove the theoretic model correct. As we know, each way of making experiments has its small flaws and inaccuracies: when we throw the coin by hand, we cannot guarantee even distribution between angular velocities, angles, etc.; robots can work incorrectly, and computer modeling algorithms are based on built-in programming languages random number generators that are sometimes not exactly random. Overall, there are different randomnesses, and using only one to prove the physical model working, would not be enough, as any given method might be biased for a specific model. But all these methods fit our theoretical model pretty well, meaning that the theoretic model is correct, and we receive the right results from it.

\section{Conclusions}
Summarising our work in facts and figures: we investigated two types of rotation – flat and volumetric and found results for a fair coin for both of them (they turned out to be different) – $$\frac{H}{R}=1.5\pm 0.05$$ and  $$\frac{H}{R}=0.8\pm 0.05$$ \\ This result differs from some of the previous findings [Hernandez-Navarro, Piñero, 2021], but also is quite similar to some others [Mahadevan, Yong, 2014]. These ratios were found for constant real parameters. We also found out that the probability of the coin landing on its sides is higher, lower k is; lower $\mu$ is. And does not depend on the mass and initial height.

Another important thing to note, perhaps as a point for future investigation, in this solution, the overall idea was to determine the end state of the coin, depending on its starting position. However, there is another way to think about it, we can try to set the starting parameters of the coin in a certain way, so it lands on a particular side. For instance, we can make the coin spin in such a way, that it is going to precess. Then it is obvious that if we can make it always face the table on one side while falling, we are going to be able to find such starting parameters, that the coin will land on a specific side.

\bibliographystyle{unsrt}  
\bibliography{references}  

Poincaré, H. (1912). Chance. The Monist, 31-52. 

Harris, T. A., and Crecelius, W. J. (1986). Rolling bearing analysis. Journal of Tribology, 108(1), 149–150. doi:10.1115/1.3261135 

Yong, E. H., and Mahadevan, L. (2011). Probability, geometry, and dynamics in the toss of a thick coin. American Journal of Physics, 79(12), 1195–1201. doi:10.1119/1.3630934

Yong, E. H., and Mahadevan, L. (2014). Statistical mechanics and shape transitions in microscopic plates. Physical Review Letters, 112(4). doi:10.1103/physrevlett.112.048101 

Stefan, R. C., and Cheche, T. O. (2016). Coin toss modeling. 
arXiv preprint arXiv:1612.06705

STRZALKO, J., GRABSKI, J., STEFANSKI, A., PERLIKOWSKI, P., and KAPITANIAK, T. (2008). Dynamics of coin tossing is predictable. Physics Reports, 469(2), 59–92. doi:10.1016/j.physrep.2008.08.003 

Kuindersma, S. R., and Blais, B. S. (2007). Teaching bayesian model comparison with the three-sided coin. The American Statistician, 61(3), 239–244. doi:10.1198/000313007x222497 

Parker, M., Rogers, J. K., and Hunt, H. (2018). Retrieved from https://www.youtube.com/watch?v=-qqPKKOU-yY 

Brown, D., and Cox, A. J. (2009). Innovative uses of video analysis. The Physics Teacher, 47(3), 145–150. doi:10.1119/1.3081296 

Keller, J. B. (1986). The probability of heads. The American Mathematical Monthly, 93(3), 191–197. doi:10.1080/00029890.1986.11971784 

Nagler, J., and Richter, P. (2008). How random is dice tossing? Physical Review E, 78(3). doi:10.1103/physreve.78.036207 

Nagler, J., and Richter, P. H. (2010). Simple model for dice loading. New Journal of Physics, 12(3), 033016. doi:10.1088/1367-2630/12/3/033016

Hernández-Navarro, L., and Piñero, J. (2021). Exact edge landing probability for the bouncing coin toss and the three-sided die problem. 
arXiv preprint arXiv:2103.10927

Mazliak, L. (2013). Poincaré and probability. Lettera Matematica, 1(1–2), 33–39. doi:10.1007/s40329-013-0004-2 
Bachelier, L. (1912). Calcul des probabilities. Paris: Gauthier-Villars. 

v. Smoluchowski, M. (1918). Über den Begriff des Zufalls und den ursprung der wahrscheinlichkeitsgesetze in der Physik. Naturwissenschaften, 6(17), 253–263. doi:10.1007/bf01491334 

Lobanov, A. I., Petrov, I. B. (2021). MATEMATICHESKOE MODELIROVANIE NELINEJNYH PROCESSOV. Publishing URAIT. [in Russian]

Reznichenko, G. Y., Myatlev, V. D., Panchenko, L. A., Terekhin, A. T. (2017). TEORIYA VEROYATNOSTEJ I MATEMATICHESKAYA STATISTIKA. MATEMATICHESKIE MODELI. Publishing URAIT. [in Russian]

Sukhov, A. D., Gaek, A. A., Coin toss modelling program in Unity with C Sharp. Version (1.0). Github. Retrieved from https://github.com/ArtemSuhov/CoinToss.

\end{document}